\documentclass{article}[12pt]
\usepackage[margin=1.25in]{geometry}

\usepackage{natbib}
 \bibpunct[, ]{(}{)}{,}{a}{}{,}%
 \def\bibfont{\small}%
\usepackage{authblk}
\usepackage[utf8]{inputenc}
\usepackage{amsmath}
\usepackage{amsthm}
\usepackage{amssymb}
\usepackage{subcaption}
\usepackage{tabularx}
\usepackage{bbm}
\usepackage{multirow}
\usepackage{graphicx}
\usepackage{caption}
\usepackage{subcaption}
\usepackage{url}

\newcommand{\X}{\mathbb{X}}
\newcommand{\pmin}{p^{\text{min}}}
\newcommand{\pmax}{p^{\text{max}}}
\newcommand{\fmin}{f^{\text{min}}}
\newcommand{\fmax}{f^{\text{max}}}
\newcommand{\ac}{\mathcal{A}}
\usepackage[boxed]{algorithm2e}
\newcommand{\p}{\mathbb{P}}

\renewcommand{\Re}{\mathbb{R}}

\usepackage{color}

\newtheorem{thm}{Theorem}

\newtheorem{defn}{Definition}

\begin{document}







\title{Learning for Constrained Optimization: \\ 
Identifying Optimal Active Constraint Sets}
\author[1]{Sidhant Misra}
\author[2]{Line Roald}
\author[3]{Yeesian Ng}
\affil[1]{sidhant@lanl.gov}
\affil[2]{roald@wisc.edu}
\affil[3]{yeesian@mit.edu}



\maketitle

\abstract{%
In many engineered systems, optimization is used for decision making at time-scales ranging from real-time operation to long-term planning. This process often involves solving similar optimization problems over and over again with slightly modified input parameters, often under tight latency requirements. 
We consider the problem of using the information available through this repeated solution process to directly learn a model of the optimal solution as a function of the input parameters, thus reducing the need to solve computationally expensive large-scale parametric programs in real time. Our proposed method is based on learning relevant sets of active constraints, from which the optimal solution can be obtained efficiently. Using active sets as features preserves information about the physics of the system, enables interpretable models, accounts for relevant safety constraints, and is easy to represent and encode.
However, the total number of active sets is also very large, as it grows exponentially with system size. 
The key contribution of this paper is a streaming algorithm that learns the relevant active sets from training samples consisting of the input parameters and the corresponding optimal solution, without any assumptions on the problem structure. 
The algorithm comes with theoretical performance guarantees, and is known to converge fast for problem instances with a small number of relevant active sets. It can thus be used to establish the practicability of the learning method.
Through extensive experiments on the Optimal Power Flow problem, we observe that often only a few active sets are relevant in practice, suggesting that the active sets is the appropriate level of abstraction for a learning algorithm to target.
}%




%



\section{Introduction}

The use of optimization methods to improve economic efficiency and enhance system performance is growing rapidly in engineering applications. Examples include the operation and planning of infrastructure like electric power grids, gas transmission systems, water and heating networks and transportation networks. Many of these systems operate in highly variable conditions, due to developments such as renewable energy integration in power systems and uncertain demand-side behaviour. To maintain safe and efficient operations, new optimized set-points are obtained by resolving the optimization problems as system conditions (and corresponding system parameters) change. The optimization problems exhibit two key features: 
\begin{enumerate}
    \item There exists a mathematical model of the physical system, including objective function and technical constraints. For optimization purposes, this model is assumed to be fully known except for those input parameters whose values are obtained in real-time, giving rise to parametric programming problems. 
    \item Constraint satisfaction corresponds to enforcing operational safety limits and adhering to physical laws, and is of paramount importance.
\end{enumerate}
As the optimization problems are solved over and over again in response to changing conditions, we build a rich history of input parameters and their corresponding optimal solutions. 
However, solving these parametric programs in real-time for realistically sized systems can be prohibitively difficult, owing to the combination of complex physical laws and tight latency requirements. This limits the frequency with which the system set-points can be updated. 


\begin{figure}[t] 
	\begin{minipage}{\textwidth}
	    \centering
			\begin{subfigure}{0.75\textwidth}
				\includegraphics[width=\linewidth]{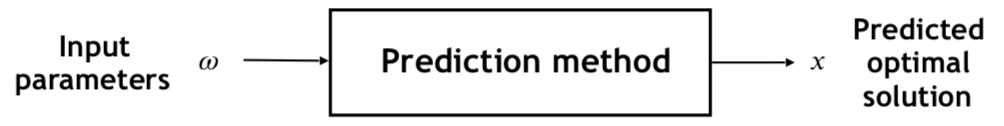}
				\caption{Predicting optimal solution from input parameters.}
				\vspace{+6pt}
                \label{fig:generic}
			\end{subfigure}
			\begin{subfigure}{\textwidth}
			    \includegraphics[width=\linewidth]{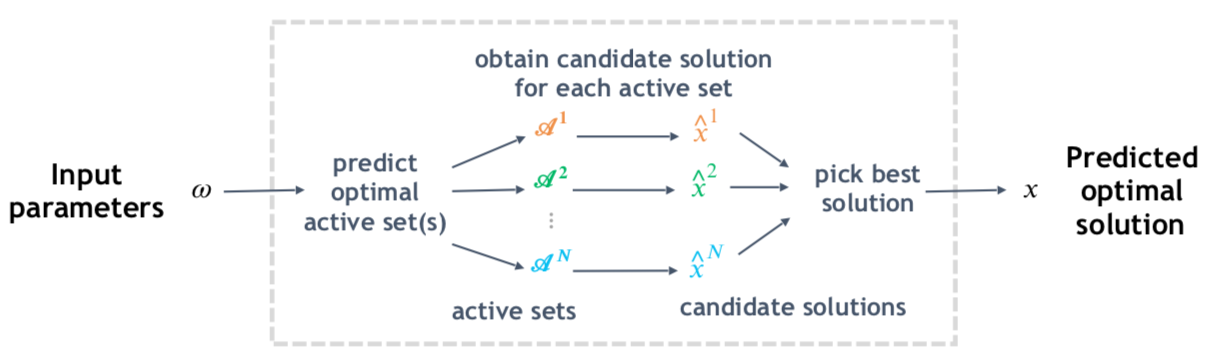}
				\caption{Predicting optimal solution using active sets as an intermediate step.}
                \label{fig:proposed}
			\end{subfigure}
	\end{minipage}
	\vspace{+4pt}
\caption{An illustration of our approach. Figure \ref{fig:generic} depicts the general approach for predicting optimal solutions from input parameters. Figure \ref{fig:proposed} depicts our approach using active sets as an intermediate step.}
\end{figure}

The goal of this paper is to develop a framework that uses machine learning methods to enable these parametric programs to be solved more efficiently. 
An intuitive approach to apply machine learning would be to use data from existing solutions, and directly learn the optimal solution from the input data, as illustrated in Figure \ref{fig:generic}. Previous literature on representing the mapping from parameters to optimal solutions in parametric programs often used reinforcement learning (RL), where the underlying physical constraints of the problem are assumed to be unknown.
We refer the reader to \cite{lillicrap2015continuous,mnih2015human,van2016deep} for examples of this approach using neural networks.
However, the neural networks mostly provide only an approximate model of the physical constraints, which may result in solutions that violate important system constraints. This has prompted a literature on ``safety problems" \citep{amodei2016concrete,garcia2015comprehensive}. Approaches to mitigate those safety concerns include modifying the exploration procedure by restricting it to a set of safe policies \citep{moldovan2012safe}, re-projecting the decisions made by the neural network into a known set of safe decisions \citep{achiam2017constrained}, or include an outer safety loop which overrides potentially unsafe decisions by the neural net \citep{fisac2017general}. Nonetheless, it remains difficult to enforce guarantees on behavioral constraints in the resulting policy \citep{moldovan2012safe, lu2017safe,achiam2017constrained} and verify that the resulting neural network policy provides a safe policy \citep{huang2017safety}. Approaches that restrict or limit the neural net also may produce sub-optimal solutions. 

A main drawback of these existing learning methods when applied in the context of engineered systems is hence their inability to enforce the constraints in the optimization accurately. One important reason for this drawback is that the methods are unable to leverage pre-existing knowledge about the mathematical form of the optimization problem. For best results, the learning methods used should exploit any known structure of the optimization problem in consideration. 

In this paper, we take a different approach to learning the optimal solution. Instead of directly learning the mapping from the input data to the optimal solution, we introduce the \emph{active set} as an intermediate step in our learning procedure. To the best of our knowledge, this has never been explored before in the literature.
The \emph{optimal active set} in an optimization problem corresponds to the set of constraints that are satisfied with equality in the corresponding optimal solution. Learning the optimal active set arguably constitute the \emph{right} level of abstraction for applying machine learning algorithm to optimization problems for several reasons:
\begin{itemize}
    \item[(i)] From an optimization perspective, the optimal active set in an optimization problem constitutes the minimal information required to recover the optimal solution. All optimization algorithms such as interior point methods, active set methods, and the simplex algorithm can see significant speed-up from the knowledge of the active set at optimality. In some cases, such as in linear programming (LP), the knowledge of the optimal active set can simplify the problem to the extent that it can be solved analytically, or in a single iteration.
    \item[(ii)]  From a machine learning perspective, using the active sets as features allows us to reduce the dimension of the learning task. Instead of learning a complex multi-dimensional, continuous-valued mapping from the input parameters to the optimal solution, the problem is converted to the simpler task of learning the mapping between the input parameters and a finite number of active sets, allowing us to exploit the known mathematical model of the system. Thus, in contrast to the standard approach in reinforcement learning, we spend no effort on learning the mathematical model itself.
    \item[(iii)] From an application perspective, active sets often correspond to meaningful operational patterns, allowing for interpretable models that can assist humans in decision-making. This point is critical in some engineering applications, where owing to the high stakes involved (such as avoiding black out in power grid), the final decision making is frequently performed by experienced system operators.
    \item[(iv)] Our intuition tells us that human engineered systems often admits just a few modes of operation. This is borne out in the numerical experiments that we ran on the optimal power flow (OPF) problem across a wide variety of transmission networks, which translates into only a small number of active sets being relevant in practice. This offers significant advantage for a machine learning approach to learn these active sets from data.
\end{itemize}


Therefore, a key contribution of this paper is to show that learning the active set is an efficient and viable approach in many practical problems, with low number of relevant active sets. To this end, we propose a statistical learning approach which provides probabilistic guarantees for out-of-sample performance, and characterizes the complexity of the task in terms of the number of active sets required to produce an optimal solution with high probability. 
The input to the algorithm is a set of training samples consisting of the input parameters and the corresponding optimal solution, without any assumptions on the problem structure. 
The algorithm is designed to terminate when a sufficient number of optimal active sets have been discovered. The algorithm also allows for the user to set termination criteria that provides theoretical performance guarantees. In the worst case, the algorithm may require an exponential number of samples. However, we prove that it will terminate fast for so-called \emph{low-complexity systems} where the number of optimal active sets is small. 
The algorithm thus simultaneously identifies the optimal active sets and establishes the practicability of the active-set-based learning method, which hinges on only a modest number of optimal active sets.
As soon as the relevant active sets are identified, we can use them to predict the optimal solution, as outlined in Figure \ref{fig:proposed}.


A related, but distinct field of study is \emph{explicit} Model Predictive Control (MPC), which was developed to improve solution times of MPC problems in an online setting with time-varying parameters. For parametric programs with a convex quadratic objective and linear constraints, the parameter space can be partitioned into a number of regions (corresponding to different optimal active sets), where the optimal solution in each region is an affine function of the input parameters.
Explicit MPC approaches generally involve an offline procedure to identify the full set of regions and corresponding optimal solutions \citep{borrelli2001efficient, Bemporad2002-cx, bemporad2002explicit}. They can then be stored in sophisticated data-structures to reduce memory storage requirements \citep{fuchs2010optimized, geyer2008optimal, jafargholi2014accelerating}. There are also advances in fast search algorithms for the ``point location problem" of identifying the corresponding active set for a given state point based on binary search trees \citep{tondel2003evaluation, johansen2003approximate, bayat2012flexible} and other data structures \citep{bayat2011using, herceg2013evaluation, zhang2016using, zhang2018kd}. 
Further, \cite{karg2018efficient} propose to design a deep neural net to specifically captures the piece-wise affine structure of the MPC solution manifold, and provides theoretical bounds on the number of neurons and layers to enable an exact representation of the solution.

An important difference to our approach is that explicit MPC searches through the full parameter space, and is hence limited to low dimensional systems, whereas our method uses information about the distribution over the input parameters to only investigate the practically relevant parts of the input space. 
Further, our approach is more general -- explicit MPC is focused on solving convex quadratic programs, whereas the type of problems under consideration in this paper belong to a more general class of problems, sometimes including non-convex or integer constraints. 


In the following we lay down the foundations for a statistical learning approach which provides probabilistic guarantees for out-of-sample performance, and characterizes the complexity of the task in terms of the number of active sets required to produce an optimal solution with high probability. 
Our main contributions in this paper are as follows:\\
(i) Establishing active sets as a natural and effective means to learn solutions to general parametric optimization problems. To the best of our knowledge, this has never been explored before in the literature.\\
(ii) Developing a streaming algorithm termed \emph{DiscoverMass} with rigorous performance guarantees to discover the collection of relevant active sets, without any \emph{a priori} assumptions on problem structure. 
The algorithm serves a dual purpose, as it both learns the relevant active sets and provides evidence of whether the suggested active-set-based learning method is applicable to a given problem instance.\\
(iii) Establishing the viability of the active-set approach on a practical problem. We employ \emph{DiscoverMass} on a set of Optimal Power Flow benchmarks representative of optimization problems solved in operation of electric transmission grids to demonstrate that the number of relevant active sets is indeed small and learning them is highly efficient in terms of the number of samples required.

We emphasize that this paper \emph{does not} provide a full end-to-end learning framework based on learning the active set. While we provide examples for how the discovered collection of active sets can be used to recover the optimal solution for a new realization of the input parameters using an ensemble policy or classification, the development of those frameworks will be considered in future work. In this paper, we lay the foundation for this future work by showing that methods based on learning the active sets are viable for practical problems. 

The remainder of the paper is organized as follows. We first describe the main idea behind learning \emph{active sets} in Section~\ref{sec:learning_for_opt}. We then develop a learning algorithm to identify a collection of active sets from training samples in Section~\ref{sec:learning-algorithm} and establish theoretical performance guarantees. We demonstrate the effectiveness of our method using an application to power systems control in  Section~\ref{sec:numerical-results}. Section~\ref{sec:conclusion} concludes the paper.

\begin{section}{Active Set Learning for Optimization} \label{sec:learning_for_opt}
The aim of this paper is to learn decision policies $x^*(\omega)$ which returns the optimal solution $x^*$ for parametric programs of the form
\begin{align} \label{eq:optimization_problem}
    x^*(\omega)\in\underset{x\in\X(\omega)}{\text{argmin}}\ f(x,\omega)
\end{align}
where $x$ denotes the decision variables, $\omega$ represent the uncertain input parameters, and
\begin{align}
    \X(\omega)=\{x\in\Re^n:\ g_j(x,\omega)\leq 0 \quad\forall j\in [m]\}
\end{align}
defines the set of feasible solutions, and $f,g_1,\dots,g_m$ are functions in $x$ for all $\omega$. We restrict ourselves to a set of parameters $\omega\in\Omega$, for which the optimization problem is actually feasible:
\begin{equation}
\Omega=\{\omega:\mathbb{X}(\omega)\neq\emptyset\}.
\end{equation}
For any given solution $x\in\mathbb{X}(\omega)$, we define the \emph{active set} corresponding to $x$ as the set of constraints $\ac$ that are binding at $x$.
For a given parameter realization $\omega$, we define the \emph{optimal active set} $\ac^*(\omega)$ as the active set at the optimal solution $x^*(\omega)$ \footnote{When there are multiple optimal solutions, the degeneracy can be handled by using any consistent tie-breaking rule,
such that the same realization $\omega$ will always be mapped to the same optimal active set.}

If we know the optimal set $\ac^*$ for a given $\omega\in\Omega$, we can solve the following optimization problem
\begin{equation} \label{eq:reduced_optimization_problem}
x^*_{\ac^*}(\omega)\in\underset{x\in\mathbb{X}_{\ac^*}(\omega)}{\text{argmin}} \ \ f(x,\omega)
\end{equation}
based on the reduced set of constraints $\ac^*$
\begin{equation} \label{eq:active_set}
\mathbb{X}_{\ac^*}(\omega) = \{x\in\Re^n:\ g_j(x,\omega) \leq 0 \quad\forall j\in \ac^*\}
\end{equation}
The solution to \eqref{eq:reduced_optimization_problem} will be feasible and optimal for the original problem \eqref{eq:optimization_problem}, but is easier to obtain since we typically have that $|\ac^*| \ll m$. 
If we have a method to map the uncertain parameters $\omega$ to the corresponding optimal active set $\ac^*$, we are thus able to obtain a the optimal solution efficiently.

Consider the case where we draw a new sample $\omega$ and do not know the optimal active set $\ac^*$, but have a collection of possible candidates $\ac_1,...,\ac_N$. Then, it is possible to obtain the optimal solution by solving the reduced problem \eqref{eq:reduced_optimization_problem} for each active set and check the resulting solutions for feasibility. Since each reduced problem \eqref{eq:reduced_optimization_problem} is a relaxation of \eqref{eq:optimization_problem}, the solutions will either be optimal or infeasible for the original problem.

\subsection{Learning the collection of relevant active sets}
While the collection of all possible active sets is finite, it is also exponential in the problem size, making it a potentially complex learning task with high sample complexity. However, in many applications such as infrastructure systems that were engineered to accommodate particular modes of operation, only a few active sets are relevant in practice. The learning task at hand thus becomes to identify those important active sets.

In our approach, the \emph{important} active sets are the active sets that have the highest probability of being optimal, under a given probability distribution over the uncertain input parameters $\omega$. Our procedure to discover the active sets is illustrated in Figure \ref{fig:learningprocedure}. We first draw i.i.d. training samples $\omega_i$ from the probability distribution over $\omega$. For each sample, we solve the optimization problem \eqref{eq:optimization_problem} to obtain the corresponding optimal solution and active set $\ac_i$, which is now marked as an \emph{observed} active set (in color in Figure \ref{fig:learningprocedure}). The collection of all active sets $\ac_1, ..., \ac_N$ that were observed during the first $\omega_i, ...,\omega_N$ samples form the set of observed active sets. The still \emph{unobserved} active sets, marked in grey in Figure \ref{fig:learningprocedure}, are the active sets which were not optimal for any of our i.i.d. samples $\omega_i$.

\subsection{From active sets to the optimal solution}
The knowledge of the relevant active sets can be used to more efficiently obtain the optimal solution. We describe two possible approaches.

\paragraph{Ensemble policy}
The ensemble policy obtains the optimal solution by solving the reduced problem \eqref{eq:reduced_optimization_problem} in parallel for each relevant active set $\ac$, and then verifying the resulting solution for feasibility in \eqref{eq:optimization_problem}. This method is described in detail in \cite{Ng2018-bj}, and illustrated in Figure \ref{fig:proposed}. 

\paragraph{Classification}
This approach is based on training a classifier to learn the mapping from the parameter $\omega$ to the corresponding optimal active set. The number of classes in this classifier is same as the number of relevant active sets. Once the optimal active set is obtained, the optimal solution can be obtained by solving the reduced optimization problem in \eqref{eq:reduced_optimization_problem}. 

Both the above approaches are simpler and more efficient when the number of potential active sets is small.
In this paper we utilize the \emph{ensemble policy} in our experiments to test the performance of our method, and leave classification as a topic for future research.

We envision any learning approach based on learning active sets to rely on combining a new realization  $\omega$ with one or more observed active sets, as illustrated in Figure \ref{fig:predictionprocedure}. 
If the optimal active set corresponds to one of the observed active sets, then the optimal solution lies within the set of optimal solutions to the corresponding collection of relaxed problems \eqref{eq:reduced_optimization_problem}. On the other hand, if the optimal active set of the new sample was previously unobserved, the solutions to the all relaxed problems \eqref{eq:reduced_optimization_problem} with $\ac_i,...,\ac_N$ may be infeasible. Since this is an undesirable outcome, we would like to bound the probability of these outcomes. 
A main contribution of this paper is to provide a learning algorithm which guarantees that the probability of encountering a previously unobserved optimal active set falls below a pre-specified level. We develop a streaming algorithm which draws new samples and observes new optimal active sets until the desired performance level can be guaranteed. This helps us distinguish between cases where we have discovered active sets which cover a smaller or larger portion of the feasible space, as illustrated in Figure \ref{fig:discovery}. The algorithm is formalized in Section~\ref{sec:learning-algorithm}.
\begin{figure}
    \centering
    \includegraphics[width=0.9\textwidth]{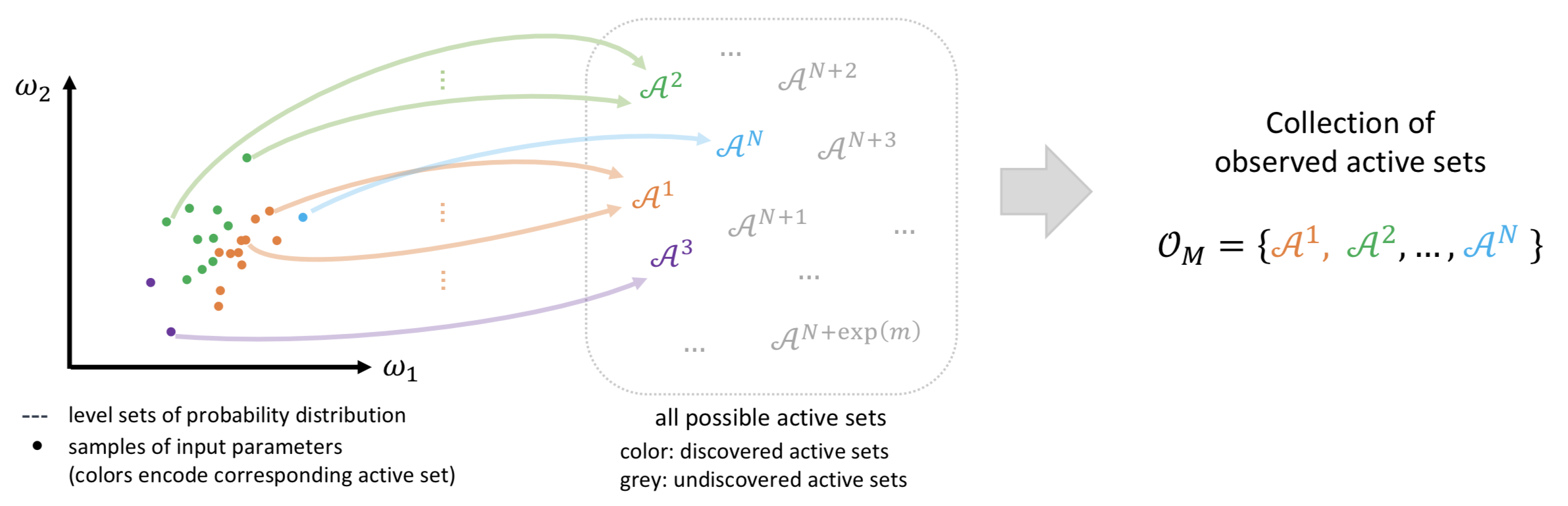}
    \caption{Learning procedure. Starting from samples (left), we discover the optimal active sets among the set of all possible active sets (middle). This provides a collection of important active sets, corresponding to those that are likely to be optimal (right).}
    \label{fig:learningprocedure}
\end{figure}
\begin{figure}
    \centering
    \includegraphics[width=0.9\textwidth]{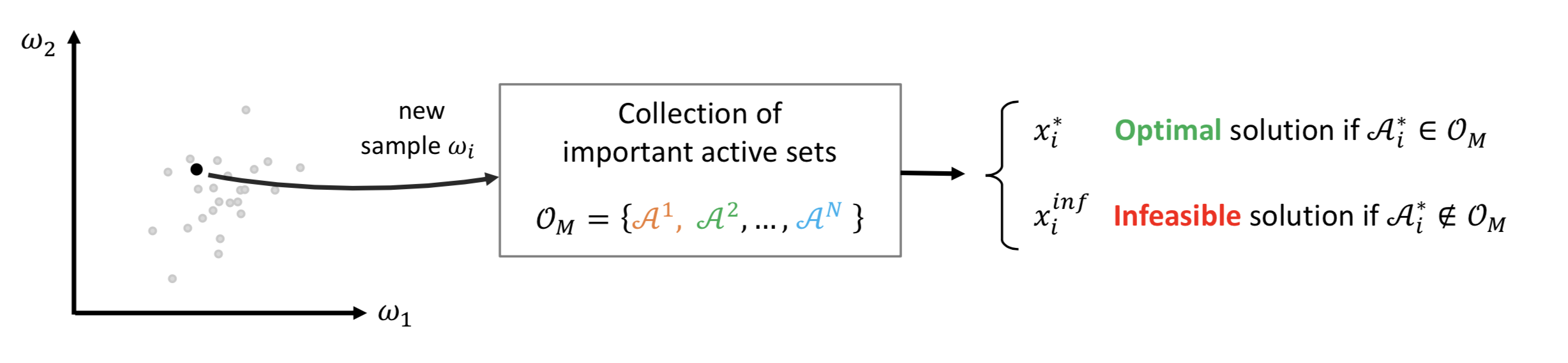}
    \caption{Prediction procedure. For a new realization $\omega$ (left), we evaluate a candidate solution for one or more observed actives sets (middle). If the optimal active set for $\omega$ corresponds to one of the observed actives sets, we recover the true optimal solution. If not, we will may end with an infeasible solution (right).}
    \label{fig:predictionprocedure}
\end{figure}

\begin{figure}[!bth] 
	    \centering
			\begin{subfigure}{0.8\textwidth}
				\includegraphics[width=\linewidth]{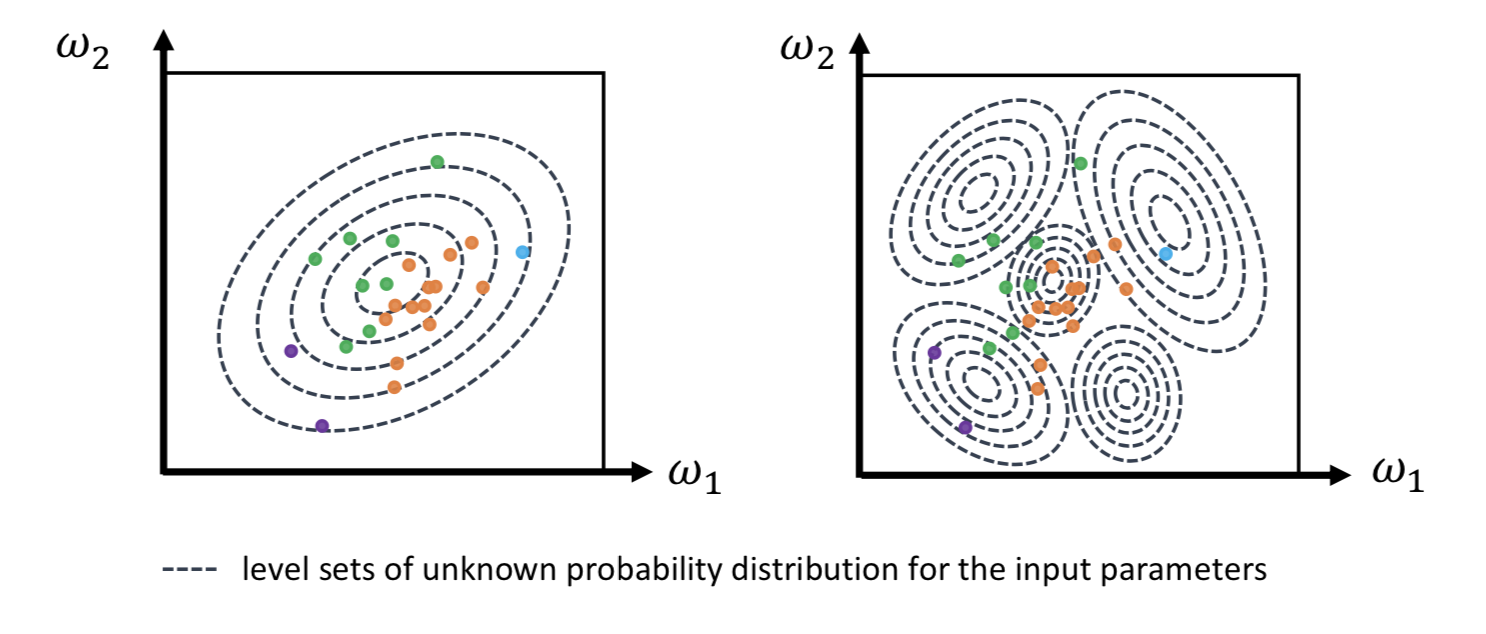}
				\caption{The samples can come from different probability distributions.}
					\vspace{+6pt}
                \label{fig:levelsets}
			\end{subfigure}

			\begin{subfigure}{0.8\textwidth}
			    \includegraphics[width=\linewidth]{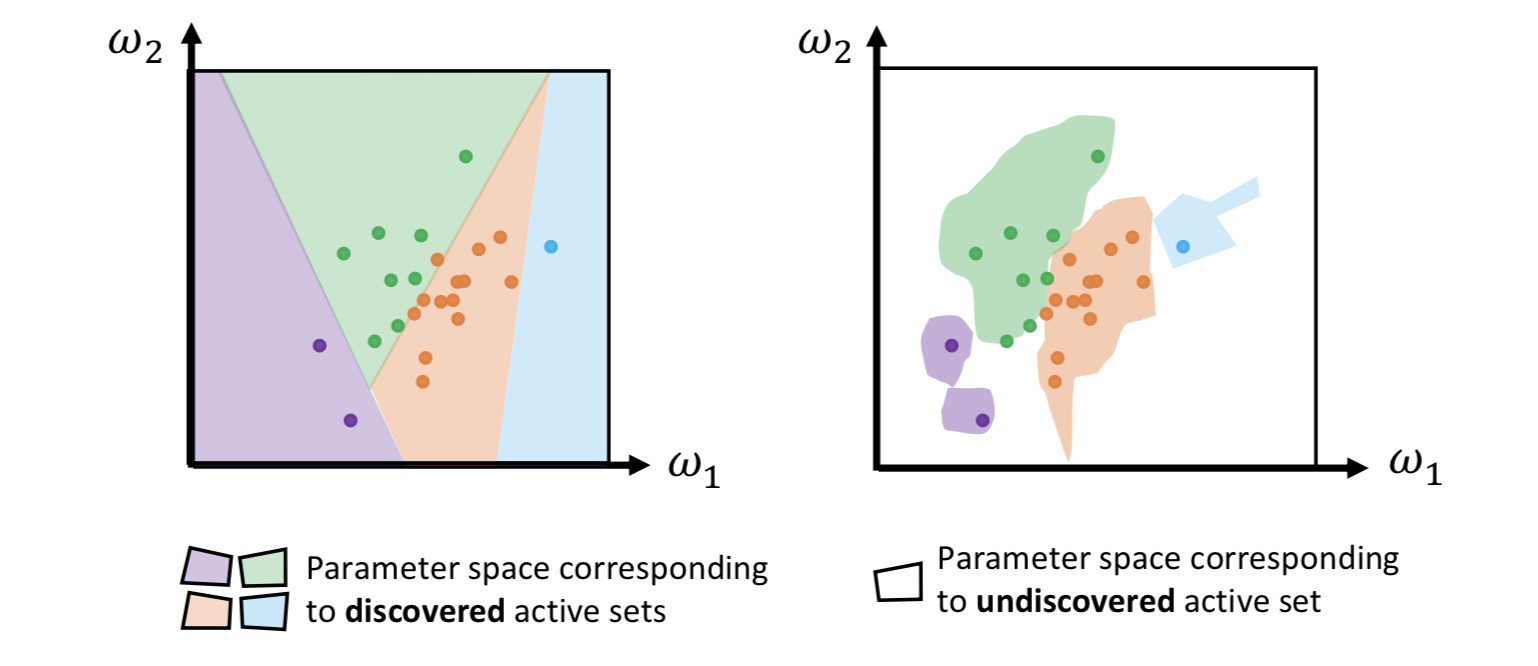}
				\caption{Observed active sets can cover a larger or smaller part of the parameter space.}
                \label{fig:activesets}
			\end{subfigure}
	\vspace{+4pt}
\caption{Our method makes no assumptions about the structure of the problem or the probability distribution we sample from, and it is hence very general. Figure \ref{fig:levelsets} shows how the samples can arise from different probability distributions. Figure \ref{fig:activesets} shows how the active sets we discovered may cover a smaller or larger region of the parameter space. In general, we do not know whether what we have observed so far corresponds to a small region (left) or a large region (right), and how much of the probability mass is included in the discovered areas. Therefore, an essential part of our learning algorithm is to identify when we have discovered ``enough'' active sets, i.e. active sets that cover a sufficient portion of the probability mass in parameter space. Note that the active sets may have an arbitrary shape in general, as we make no assumptions on the problem structure. \label{fig:discovery}}
\end{figure}

\end{section}

\begin{section}{Algorithm to learn relevant active sets} \label{sec:learning-algorithm}
In this section, we describe a streaming algorithm that enables us to learn a collection of active sets whose combined probability mass exceeds a user-defined safety level. 
The algorithm draws random samples of the uncertain input parameters $\omega_i$ from the distribution $\p_{\omega}(.)$, solves the optimization problem \eqref{eq:optimization_problem} to obtain the corresponding active set $\ac_i$ and uses a \emph{stopping criterion} to ensure that sufficiently many active sets have been discovered. We provide details of the algorithm, and theoretical guarantees on its performance below.

Let $\mathcal{B}$ denote the finite set of all possible active sets, which is exponential in the number of constraints in the problem. 
We define \emph{observed} and \emph{unobserved} active sets in the following way. 
\begin{defn} \label{def:unobserved_set}
Let $\omega_1, \ldots, \omega_M$ be $M$ i.i.d. samples drawn from the uncertainty distribution and $\ac_1, \ldots, \ac_M$ denote the optimal active sets corresponding to each $\omega_i$. We call $\mathcal{O}_M = \cup_{i=1}^M \{\ac_i\}$ the set of observed active sets, and $\mathcal{U}_M = \mathcal{B} \setminus \mathcal{O}_M$ the set of unobserved active sets.
\end{defn}
In Figure~\ref{fig:activesets}, $\mathcal{O}_M$ and $\mathcal{U}_M$ correspond to the colored and white regions respectively. The streaming algorithm iteratively increases the number of samples $M$ until the observed set $\mathcal{O}_M$ contains most of the important active sets. A crucial aspect of the algorithm is that knowledge of neither the shape and number of the active sets (colored regions) nor the number of samples $M$ is required a-priori, but is decided on-the-fly. The algorithm therefore is guaranteed to converge for any system if enough samples are available, but can also be stopped at any time before termination with a lower safety guarantee.

In the following, we use the notion that the {importance} of an active set or a collection of active sets corresponds to the \emph{mass} contained in it.
\begin{defn} \label{def:coverage}
The \emph{mass} of an active set is defined as the probability that it is optimal for \eqref{eq:optimization_problem} under the distribution $\p_{\omega}(.)$ on $\omega$. For $\ac \in \mathcal{B}$,
\begin{align}
    \pi(\ac) = \p_{\omega}(\ac = \ac^*(\omega)) 
\end{align}
For any subset $\mathcal{S} \subseteq \mathcal{B}$, we define the \emph{mass} of $\mathcal{S}$ as $\pi(\mathcal{S}) = \sum_{\ac \in \mathcal{S}} \pi(\ac)$.
\end{defn}
After drawing $M$ samples and observing a collection of active sets $\mathcal{O}_M$, we compute the rate of discovery. The rate of discovery quantifies the fraction of samples that correspond to observing an active set that was not observed among the first $M$ samples, and is formally defined below.
\begin{defn} \label{def:rod}
Let $W$ be a positive integer denoting the window size. Let $\omega_1, \ldots, \omega_{M+W}$ be $M+W$ i.i.d. samples drawn from the uncertainty distribution and let $\ac_i$ denote the optimal active set for $\omega_i$. We denote by $X_i$ the random variable that encodes whether a ``new active set", i.e., an active set not observed within the first $M$ samples, was observed in the $(M+i)^{th}$ sample, i.e.
 \begin{align}
        X_i = \begin{cases}
                    1, \ \mbox{if} \ \ac_{M+i} \notin \{\ac_1\} \cup \ldots \cup \{\ac_{M}\}, \\
                    0, \ \mbox{otherwise}.
               \end{cases}
    \end{align}
    Then the \emph{rate of discovery} over the window of size $W$ is given by $\mathcal{R}_{M,W}$ and is defined as
    \begin{align} \label{eq:R_def}
        \mathcal{R}_{M,W}  = \frac{1}{W} \sum_{i=1}^{W} X_i.
    \end{align}
\end{defn}

The rate of discovery $\mathcal{R}_{M,W}$ is directly related to the mass of the undiscovered set $\mathcal{U}_M$ of active sets. Intuitively, if we have already discovered the active sets that contain most of the mass, then the rate of discovery is expected to be low since most of the time we are likely to observe one of the active sets we have already observed so far. The following theorem makes this notion rigorous.
\begin{thm} \label{thm:rod_proximity}
Let $\{\omega_i\}_{i=1}^{\infty}$ be a sequence of i.i.d. observations and define the unobserved set and rate of discovery as in Definition~\ref{def:unobserved_set} and \ref{def:rod} respectively. Then the following statements hold.
Let the probability mass in the unobserved set be denoted by $\pi(\mathcal{U}_M)$ and let the window size be such that 
    \begin{align} \label{eq:window_size}
        W = W_M \geq c \max \{\log \underline{M},  \log M\}.
    \end{align}
    Then we can bound the probability mass in the undiscovered set by
    \begin{align}\label{eq:confidencedelta}
        \mathbb{P}& \left(\pi(\mathcal{U}_M) - \mathcal{R}_{M,W}  \leq  \epsilon \quad \forall M \geq 1 \right) > 1 - \frac{ \gamma}{\gamma-1} \frac{1}{(\underline{M}-1)^{\gamma-1}}, \quad \mbox{where } \gamma = \frac{c \epsilon^2}{2}.
    \end{align}
\end{thm}
The proof is deferred to Section~\ref{sec:proofs}.
Theorem~\ref{thm:rod_proximity} forms the basis of the stopping criterion of our streaming algorithm. In the following subsection, we present this algorithm and provide rigorous guarantees on its performance based on the rate of discovery criterion.

\begin{subsection}{Streaming Algorithm for Learning Important Active Sets} \label{subsec:algorithm}
We present the following algorithm called \emph{DiscoverMass} which takes in as input parameters the quantities $\alpha, \epsilon, \delta >0$. Here, $\alpha$ represents the limit on the undiscovered mass and requires that the algorithm returns a set of active sets of mass at least $\alpha$. The quantity $\epsilon$ represents the difference between the undiscovered mass and the rate of discovery, and $1-\delta$ represents the confidence.
We can also choose the hyperparameter $\gamma$. We have explicitly left the dependence on all the parameters in the algorithm description so that interested readers can choose the values that may be suitable for their application. We specify the concrete set of chosen values along with a brief analysis of the effect of tuning the hyperparameters in the numerical experiment section.
 \begin{algorithm}[H]
        \SetAlgoLined
        \KwData{$\alpha, \epsilon, \delta, \gamma$}
        Compute $c = 2\gamma/\epsilon^2$ and $\underline{M} = 1 + \left(\frac{\gamma}{\delta(\gamma-1)}\right)^{\frac{1}{\gamma-1}}$\; 
        Initialize $M = 1$ and $\mathcal{O} = \emptyset$\; 
        \Repeat{$\mathcal{R}_{M,W_M} < \alpha - \epsilon$}{
        Calculate window size $W_M= c \max \{\log \underline{M}, \log M \}$ \;
        Draw $1  + W_M - W_{M-1}$ additional samples from $\mathbb{P}_{\omega}$ to obtain a total of $M + W_M$ samples \; 
        Solve optimization problem \eqref{eq:optimization_problem} for each new sample\;
        Add the newly observed active sets to $\mathcal{O}$ \; 
        Compute $\mathcal{R}_{M,W_M}$ as given in \eqref{eq:R_def}\; 
        Update M = M+1\;
        }
        \Return $\mathcal{O},M, R_{M,W_M}$
    \caption{DiscoverMass}
    \end{algorithm}
    
\emph{DiscoverMass} terminates with a set of observed active sets $\mathcal{O}$ that has at least $1-\alpha$ of the mass with high probability. The following theorem guarantees that with the stopping criterion described above, \emph{DiscoverMass} indeed succeeds in meeting the requirements imposed by the input parameters. 

\begin{thm} \label{thm:stopping_criterion}
If the algorithm \emph{DiscoverMass} is employed with input parameters $\alpha, \epsilon, \delta, \gamma$, then with probability at least $\mathbf{1-\delta}$ the algorithm will terminate having discovered active sets with mass at least $\pi(\mathcal{O}) \geq 1 - \alpha$.
\end{thm}
The proof of Theorem \ref{thm:stopping_criterion} can be found in Section \ref{sec:proofs}.
A key advantage of \emph{DiscoverMass} supported by Theorem~\ref{thm:stopping_criterion} is that it makes absolutely no \emph{a priori} assumption regarding the underlying distribution of active sets in the system. This makes \emph{DiscoverMass} a rigorous tool to quantify the unknown underlying distribution.

    
\end{subsection}

\begin{subsection}{Low-complexity systems} \label{subsec:low-complexity_systems}
Although Theorem~\ref{thm:stopping_criterion} guarantees that the algorithm \emph{DiscoverMass} will almost always succeed in finding the desired amount of mass, it gives no guarantees on the number of steps $M$ the algorithm needs to discover this mass. In fact, the number of steps is highly dependent on the properties of the system and can vary significantly across systems. For example in a system where there are $2^N$ active sets all of which are equally probable, i.e. $\p_{\omega}(\ac_i) = 2^{-N}$ for all $i \in [2^N]$, then it is clear that discovering $(1-\alpha)$ fraction of the mass requires at least $M > (1-\alpha)2^N$ samples. However, for some systems, which we refer to as low-complexity systems, it is possible to discover a large fraction of the mass quickly. Fortunately, many practical systems fall within this category. In Section~\ref{sec:numerical-results}, we have found that many practical systems can be described as \emph{low-complexity} as it is possible to discover a large fraction of the mass quickly. We quantify this intuition with the following definition.
    
    \begin{defn} \label{def:low-complexity_systems}
        We say that a system is a \textbf{low-complexity system} with parameters $\alpha_0$ and $K_0$ if there are at most $K_0$ active sets that contain at least $1-\alpha_0$ fraction of the mass. More precisely, there exists active sets $\ac_{i_1}, \ldots, \ac_{i_{K_0}}$ such that $\sum_{j=1}^{K_0} \pi\left(\ac_{i_j}\right) > 1-\alpha_0$.
    \end{defn}
    
    The next theorem guarantees that for \emph{low-complexity systems} the algorithm terminates within a few iterations.
    \begin{thm} \label{thm:low-complexity_systems}
        Suppose that we run the algorithm \emph{DiscoverMass$(\alpha,\epsilon,\delta,\gamma)$} on a low-complexity system with parameters $(\alpha_0,K_0)$ with $\alpha >  \alpha_0$. Then with probability at least $1 - \delta - \delta_0$ the algorithm terminates in less iterations than $M = \frac{1}{\alpha - \alpha_0}\left(K_0 \log 2 + \log 1/\delta_0 \right)$.
    \end{thm}
    The proof of Theorem \ref{thm:low-complexity_systems} is found in the next section. From Theorem~\ref{thm:low-complexity_systems}, we see that the sample complexity of discovering active sets with sufficient mass for low-complexity systems scales at most linearly with the number of underlying important active sets $K_0$. 
    
  \end{subsection}
\end{section}

\begin{section}{Proofs} \label{sec:proofs}
In this section, we provide the proofs of the three theorems stated in Section~\ref{sec:learning-algorithm}.
\subsection{Proof of Theorem~\ref{thm:rod_proximity}}
We first prove part $(a)$.  For any $M$ we can bound the probability that the rate of discovery is far from the mass of the unobserved set by 
 \begin{align}
            \mathbb{P} & (\pi(\mathcal{U}_M) - \mathcal{R}_{M,W}  >  \epsilon ) \nonumber \\
            &= \sum_{u} \mathbb{P}  (\pi(\mathcal{U}_M) - \mathcal{R}_{M,W}  >  \epsilon \mid \mathcal{U}_M = u)  \mathbb{P}(\mathcal{U}_M = u) \nonumber \\
            &= \sum_{u} \mathbb{P}  \left(  \pi(u) -  \frac{\sum_{i=1}^{W_M} X_i }{W_M}>  \epsilon \mid \mathcal{U}_M = u\right)  \mathbb{P}(\mathcal{U}_M = u), \label{eq:partial}
        \end{align}
        where the random variables $X_i$ are defined in \eqref{eq:R_def} and conditioned on $\mathcal{U}_M = u$ are i.i.d. Bernoulli random variables with mean $\pi(u)$. Using standard large deviation inequalities, we can bound the probability in \eqref{eq:partial} as
        \begin{align}
            \mathbb{P}  (\pi(\mathcal{U}_M) - \mathcal{R}_{M,W}  >  \epsilon )  &< \sum_{u} e^{-W_M\epsilon^2/2} \mathbb{P}(\mathcal{U}_M = u)  \nonumber \\
            &= e^{-W_M\epsilon^2/2}. \label{eq:fixed_M_bound}
        \end{align}
        Therefore,
        \begin{align}
            \mathbb{P}& (\exists \ M  \geq 1  \quad  \mbox{s.t.} \quad  \pi(\mathcal{U}_M) - \mathcal{R}_{M,W}  >  \epsilon ) \nonumber \\
            & \leq \sum_{M=1}^{\infty} \mathbb{P} (\pi(\mathcal{U}_M) - \mathcal{R}_{M,W}  >  \epsilon ) \nonumber \\
            &= \sum_{M=1}^{\underline{M}} \mathbb{P} (\pi(\mathcal{U}_M) - \mathcal{R}_{M,W}  >  \epsilon ) \nonumber \\
            &+  \sum_{M=\underline{M}}^{\infty} \mathbb{P} (\pi(\mathcal{U}_M) - \mathcal{R}_{M,W}  >  \epsilon )  \label{eq:partial_sums} \\
             & < \sum_{M=1}^{\underline{M}}  e^{-c \log \underline{M}\epsilon^2/2} + \sum_{M=\underline{M}}^{\infty} e^{-c \log M \epsilon^2/2} \\
            & \leq \frac{1}{\underline{M}^{\frac{c \epsilon^2}{2} -1}} + \sum_{M=\underline{M}}^{\infty}\frac{1}{\underline{M}^{\frac{c \epsilon^2}{2}}} \leq  \frac{\gamma}{\gamma-1} \frac{1}{(\underline{M}-1)^{\gamma-1}}.\nonumber
        \end{align}
       
           
   
\subsection{Proof of Theorem~\ref{thm:stopping_criterion}}
    Let $\mathbb{E}$ be the event that the algorithm terminates after having discovered smaller than $\epsilon$ fraction of the mass. Then
        \begin{align}
            \mathbb{P}\left( \mathbb{E} \right) &\leq \mathbb{P}\left(\exists M \geq 1 \mid \mathcal{R}_{M,W} < \epsilon - \delta, \pi(\mathcal{U}_M)  > \epsilon \right) \\
            & \leq \mathbb{P} \left(\exists M \geq 1 \mid \pi(\mathcal{U}_M) - \mathcal{R}_{M,W}  > \epsilon  \right) \\
            & \stackrel{(a)}{<} \frac{\gamma}{\gamma-1} \frac{1}{(\underline{M}-1)^{\gamma-1}} \stackrel{(b)}{=} \delta,
        \end{align}
        where the last inequality $(a)$ follows from Theorem~\ref{thm:rod_proximity}, and $(b)$ follows from the value of $\underline{M}$ in the initialization step in \emph{DiscoverMass}.
        
\subsection{Proof of Theorem~\ref{thm:low-complexity_systems}}
	    Let $I^* = \{\ac_{i_1}, \ldots, \ac_{i_{K_0}} \}$ denote the set of $K_0$ active sets defined in Definition~\ref{def:low-complexity_systems} that contain at least $1-\alpha_0$ of the mass.  For any $M \geq 1$ recall that $\mathcal{U}_M$ denotes the set of unobserved active sets as in Definition~\ref{def:unobserved_set}. We define the sets  $\mathcal{U}_{M_{I^*}} = \mathcal{U}_M \cap I^*$ and $\mathcal{U}_{M_{\bar{I}^*}} = \mathcal{U}_M \cap \bar{I}^*$ where $\bar{I}^*$ denotes the complement of the ensemble $I^*$, i.e., the set of active sets that are not a part of $I^*$. It follows that
	    $\mathcal{U}_M = \mathcal{U}_{M_{I^*}} \cup \mathcal{U}_{M_{\bar{I}^*}}$ and $\pi(\mathcal{U}_M) = \pi(\mathcal{U}_{M_{I^*}})  + \pi(\mathcal{U}_{M_{\bar{I}^*}}) $.
	    We compute
	    \begin{align}
	        \mathbb{P}&\left(\pi(\mathcal{\mathcal{O}_M}) < 1 - \alpha \right) = \mathbb{P}(\pi(\mathcal{U}_M) > \alpha) \nonumber \\
	        & = \mathbb{P}(\pi(\mathcal{U}_{M_{I^*}})  + \pi(\mathcal{U}_{M_{\bar{I}^*}}) > \alpha) \nonumber \\
	        & \leq \mathbb{P}(\pi(\mathcal{U}_{M_{I^*}}) > \alpha-\alpha_0) \nonumber \\
	        &= \sum_{u \subseteq I^*} \mathbb{P}(\mathcal{U}_{M_{I^*}} = u)   \ \mathbbm{1}_{\{\pi(u) > \alpha-\alpha_0\}} \nonumber \\
	        &\leq  \sum_{u \in I^*} \left(1- (\alpha-\alpha_0) \right)^M \nonumber \\
	        & \leq 2^K e^{-(\alpha-\alpha_0)M}  \stackrel{(a)}{<} \delta_0, \label{eq:low-complexity_intermediate}
	    \end{align}
	    where $(a)$ follows from the choice of $M$ in the theorem. The proof follows by combining \eqref{eq:low-complexity_intermediate} with Theorem~\ref{thm:rod_proximity}.
\end{section}

\begin{section}{Numerical results} \label{sec:numerical-results}

We test our learning framework on the DC Optimal Power Flow (OPF) problem, an optimization problem widely used in electricity market clearing as well as operation and planning of electric transmission grids. 
For real time operation, the reaction of generators to fluctuations in renewable energy must be calculated reliably within a short period of time, so it is an example with high practical relevance that is subject to stringent limits on computational time. 
A simple version of the DC OPF can be stated as follows
%
\begin{subequations} \label{eq:opf}
	\begin{align}
        \rho(\omega)\in\underset{p\in\Re^n}{\text{argmin}}\ &c^\top p \\
	    \text{s.t.}\ &e^\top(p-d+\omega) = 0 \label{eq:powerbal}\\
    	&\pmin\leq p\leq\pmax \label{eq:genlim}\\
    	&\fmin\leq M(Hp-d+\omega)\leq\fmax  \label{eq:linelim}
	\end{align}
\end{subequations}
Here, decisions $p\in\Re^n$ are made on the active power generation at each generator. The aim is to provide power at minimum generation cost $c^\top p$ to a set of loads $d\in\Re^{v}$ with random real-time variations $\omega$. 
The solution is subject to a total power balance constraint \eqref{eq:powerbal}, as well as limits both on the minimum and maximum power generation $\pmin,\pmax\in\Re^n$ at each generator \eqref{eq:genlim}, and the minimum and maximum admissible flow $\fmin,\fmax\in\Re^m$ on transmission lines \eqref{eq:linelim}. The topology of the network is encoded in the matrix $H\in\Re^{v\times n}$ mapping the power from each generator to their corresponding bus, and the matrix of power transfer distribution factors $M\in\Re^{m\times v}$ \citep{christie2000transmission}. Here $e$ is the vector of ones.


\begin{subsection}{Test set-up}
We test our algorithm by running extensive simulations across a variety of networks \citep{li2010small,lesieutre2011examining,grigg1999ieee,birchfield2017grid,price2011reduced} from the IEEE PES PGLib-OPF v17.08 benchmark library \citep{pglib_opf}. We report general results for 15 different test cases, varying in size from 3 to 1951 buses.\\
For each system, we consider two different distributions for the uncertain load deviations $\omega$.
First, we assume a \emph{multivariate normal distribution}, where $\omega$ is defined as random vector of independent, zero mean variables with standard deviations $\sigma=0.03d$ and zero correlation between loads. 
Second, we assume that each entry $\omega_i$ follows a \emph{uniform distribution} with support $\omega_i\in[-3\sigma, +3\sigma]=[-0.09d,0.09d]$. The uniform distribution is intentionally designed to spread the probability mass more evenly across a larger region, thus enabling us to compare cases with a different number of relevant active sets and different probability mass in each active set. Note that our method is able to handle any distribution, as long as a sufficient number of training samples is available.\\
The algorithm \emph{DiscoverMass} has four input parameters, $\alpha,~\delta,~\epsilon,~\gamma$. We set the maximum mass of the undiscovered bases $\alpha=0.05$ and the confidence level $\delta=0.01$. 
The remaining parameters are set to $\epsilon=0.04$ and $\gamma=2$.
Note that $\epsilon$ and $\gamma$ are hyperparameters that can be tuned, leading to different minimum window size and $\underline{M}$. For details about the choice of parameters and their effects, we refer the reader to supplementary material. 
We run the algorithm until termination, or until $M=22'000$. 
%
The procedure is implemented in Julia v0.6 \citep{julia}, using JuMP v0.17 \citep{DunningHuchetteLubin2017}, PowerModels.jl v0.5 \citep{power_models} and OPFRecourse.jl \citep{Ng2018-bj}.
\end{subsection}

	\begin{table*}[t!]
    \centering
    \scriptsize
    \renewcommand{\arraystretch}{1.2}
    \begin{tabular}{lccccc||ccccc}
    \multicolumn{1}{c}{} & \multicolumn{5}{c}{\footnotesize\textbf{Normal distribution}} &  \multicolumn{5}{c}{\footnotesize\textbf{Uniform distribution}} \\ 
    \hline
                        & $K_M$   & $M$  & $W_M$   & $R_{M,W}$ & $\mathbb{P}(p^*)$     & $K_M$ &  $M$ & $W_M$  & $R_{M,W}$ & $\mathbb{P}(p^*)$\\
    \hline
    &&&&&&&&&&\\[-4pt]
    \textbf{Low-Complexity} &&&&&&&&&&\\
    \hline
    case3\_lmbd         & 1     & 1      & 13'259  & 0.0     & 1.0       & 1     & 1      & 13'259  & 0.0    & 1.0 \\
    case5\_pjm          & 1     & 1      & 13'259  & 0.0     & 1.0       & 1     & 1      & 13'259  & 0.0    & 1.0\\
    case14\_ieee        & 1     & 1      & 13'259  & 0.0     & 1.0       & 1     & 1      & 13'259  & 0.0    & 1.0\\
    case30\_ieee        & 1     & 1      & 13'259  & 0.0     & 1.0       & 1     & 1      & 13'259  & 0.0    & 1.0  \\
    case39\_epri        & 2     & 2      & 13'259  & 0.0     & 1.0       & 2     & 2      & 13'259  & 0.0008 & 0.9998   \\
    case118\_ieee       & 2     & 33     & 13'259  & 0.0     & 1.0       & 2     & 4      & 13'259  & 0.0019 & 0.9984   \\
    case57\_ieee        & 2     & 2      & 13'259  & 0.0003  & 0.9997    & 3     & 46     & 13'259  & 0.0    & 1.0\\
    case1888\_rte       & 3     & 6      & 13'259  & 0.0     & 1.0       & 3     & 10     & 13'259  & 0.0    & 1.0 \\
    case1951\_rte       & 5     & 47     & 13'259  & 0.0069  & 0.9943    & 11    & 63     & 13'259  & 0.0084 & 0.9901    \\
    case162\_ieee\_dtc  & 7     & 91     & 13'259  & 0.0054  & 0.9925    & 17    & 192    & 13'259  & 0.0085 & 0.9926    \\
    case24\_ieee\_rts   & 10    & 1456   & 18'209  & 0.0     & 1.0       & 11    & 64     & 13'259  & 0.0047 & 0.9941  \\[+2pt]
    \textbf{High-Complexity} &&&&&&&&&&\\
    \hline
    case73\_ieee\_rts   & 19    & 1258   & 17'844  & 0.0087  & 0.9931    & 130   & 22'000 & 24'977 & 0.0136  & -  \\
    case300\_ieee       & 24    & 1257   & 17'842  & 0.0073  & 0.9919    & 293   & 9095   & 22'789  & 0.0099 & 0.9897    \\
    case200\_pserc      & 174   & 4649   & 21'112  & 0.0099  & 0.9909    & 236   & 6741   & 22'040  & 0.0099 & 0.9901    \\
    case240\_pserc      & 2993  & 22'000 & 24'997  & 0.0795  & -         & 2993  & 22'000 & 24'997  & 0.0795 & -   \\
    \end{tabular}
    \caption{Outcome of the learning algorithm for different systems for normal and uniformly distributed parameters, sorted by the number of discovered active sets. The number in the case name indicates the number of buses in the system. $K_M$: Number of active sets. $M$: Number of samples until termination. $W_M$: Window size. $R_{M,W}$: Rate of discovery. $\mathbb{P}(p^*)$: Probability of obtaining the optimal solution.}
    \label{tab:discoverMass_results}
    \end{table*}
    
\subsection{Parameter choices in experiments}
Varying the hyper-parameter $\gamma$ results in different values for the window size $W_M$ and minimum $\underline{M}$, while changing $\epsilon$ influences both the stopping criterion $R_{M,W}$ and the size of the window $W$. We note that there is a trade-off in the choice of both of those parameters. If $\gamma$ is chosen such that $\underline{M}$ is low, the initial window size is smaller, but will rapidly start increasing as soon as $M>\underline{M}$. For $\epsilon$ a lower value implies a steep increase in the window size, but also reduces the requirements on the rate of discovery.

The difference between the probability mass of the undiscovered bases and the rate of discovery $R_{M,W}$ in the window $W$ are set to $\epsilon = 0.04$, while we choose the hyperparameter $\gamma=2$. Together, this implies the stopping criterion $R_{M,W}\leq0.01$. The minimum window size $W_{M}\geq 13'259$ samples, and remains constant as long as $M\leq\underline{M}=201$ samples. We run the algorithm until termination, or until $M=22'000$. To assess how accurately we capture the probability of the undiscovered active sets, we run an out-of-sample test with 20'000 samples. 

\begin{subsection}{Numerical results}
Table~\ref{tab:discoverMass_results} lists the simulation results for the 15 different networks with sizes ranging from 3 to 1951 nodes. For each network, we report the number of samples to termination $M$, the number of active sets $K_M$ discovered within the first $M$ samples and the window size $W_M$ required to guarantee a confidence of $\alpha<0.05$. 
In addition, we list the rate of discovery $R_{M,W}$ at termination. All results are included for both the normal and uniform distributed parameters, and the system is sorted based on the number of active sets. 
Note that the choice of input parameters leads to a minimum window size $W_{M}\geq 13'259$ samples, and remains constant as long as $M\leq\underline{M}=201$ samples. 

\paragraph{Low-complexity systems} We observe that for most systems, the learning algorithm performs well. Most systems have a relatively low number of relevant active sets ($<10$), and terminate after less than 200 samples while the window size is still $W_M = W_{M,min} = 13'259$. System size is not a determining factor for the number of relevant active sets, with even the largest systems case1888\_rte and case1951\_rte exhibiting a low number of relevant active sets.
This validates the intuition that many engineered systems (at least arising in the power systems benchmarks we've looked at) can be considered \emph{low-complexity} as defined in Section~\ref{subsec:low-complexity_systems}. Comparing $M$ and $K_M$ in Table~\ref{tab:discoverMass_results}, we see a strong correlation between the number of active sets with majority of the mass and the samples till termination $M$. This is in agreement with the sample complexity derived in Theorem~\ref{thm:low-complexity_systems}. The fact that most systems are low-complexity suggests that our proposed method for parameterizing and learning an optimal policy based on active sets might be effective in a wide variety of settings arising in engineered systems.\\
We observe that the primary factor affecting the number of samples till termination is not the number of important active sets, but the existence of active sets with small but not insignificant probability mass (typically in the range 0.001-0.01). Due to their low probability mass, those active sets require a large number of samples for discovery, while the cumulative probability of the still undiscovered active sets is still high enough to exceed the stopping criterion $R_{M,W}\leq0.01$.

\paragraph{High-complexity systems} The outliers in this analysis are cases such as case24\_ieee\_rts, case73\_ieee\_rts, case200\_pserc and case240\_pserc.
For case73\_ieee\_rts (uniform) and case240\_pserc (normal and uniform), the algorithm does not terminate until we reach the upper limit of $M=22'000$. For case240\_pserc, the rate of discovery remains above $R_{M,W}\geq0.07$ even after the discovery of 2993 active sets. 
We believe that this performance deficiency is due to atypical system characteristics. The two systems case200\_pserc and case240\_pserc, are the result of power grid models that have been subjected to network reduction, hence altering the properties of the systems relative to the other, more typical systems. Similarly, case24\_ieee\_rts and case73\_ieee\_rts are special in that they have a large number of generators relative to the number of nodes, and case73\_ieee\_rts consists of three identical areas. Due to the nature of their construction, these instances are useful adversarial systems for analyzing the performance of different learning algorithms.


%

\paragraph{Influence of parameter distribution} We observe that the uniform distribution typically discovers a larger number of active sets compared to the normal distribution and requires more samples $M$ before termination. 
With the probability mass spread more evenly across the probability space, there are a larger number of active sets with significant probability mass, and fewer active sets with high probability mass (for the multivariate normal, there is a concentration of mass in active sets closer to the forecasted value). However, for systems like case24\_ieee\_rts, the uniform distribution appears to encourage earlier termination due to a faster exploration rate. 
\paragraph{Out of sample prediction performance}
The learning algorithm is compared against a full optimization approach where for each realization of $\omega$ the problem \eqref{eq:optimization_problem} is solved using standard optimization methods. This approach has higher computational requirements, but produces optimal solutions with probability one. Our learning based method has better computational efficiency but has a non-zero probability of failure, corresponding to returning infeasible solutions.
To assess the performance of our learning algorithm we evaluated the predictions generated by the ensemble policy (where the solution is evaluated for all discovered active sets, then checked for feasibility) on $20'000$ test samples. \\
The results are shown in Table~\ref{tab:discoverMass_results} which shows the probability of failure. We see that this value is consistently lower than the mass of the undiscovered set $\alpha$ provided as input to the \emph{DiscoverMass} algorithm.
Intuitively, this is to be expected, since Theorem~\ref{thm:rod_proximity} strives to provide  upper bounds on the mass of the undiscovered set that are universal in nature and are guaranteed to work for all almost any continuous optimization problem and various kinds of uncertainty, without pre-specifying the number of samples $M$. 
\end{subsection}

\end{section}

\begin{section}{Conclusion}
\label{sec:conclusion}
We present an approach to learn the mapping from a set of uncertain input parameters to the optimal solution for parametric programs. Our target application is engineered systems, where optimization problems used in operational decision making are typically solved repeatedly, but with varying input parameters. Our approach is based on learning the active sets of constraints at the optimal solution. By relying on active sets, we are able to decompose the learning task into a simpler classification task, and at the same time exploiting the knowledge of the exact mathematical model of the optimization problem, and enabling tight constraint enforcement. The viability of learning active sets and the performance of our approach is validated using a series of experiments based on benchmark optimization problem used in electric grid operation.
The prediction step in this paper used the so-called ensemble policy. In future work we will explore classification and ranking algorithms to learn the mapping between the realization of the uncertain parameter onto the optimal active set. Further topics of exploration include application to long-term planning and multi-stage optimization problems. 
\end{section}

\bibliographystyle{informs2014}
\bibliography{references}

\end{document}